\documentclass[a4paper,11pt]{article}
\usepackage{amsmath,amssymb,amsthm,amsfonts,bbm}
\usepackage{graphics,psfrag,epic,eepic}

\title{Random recursive trees and the Bolthausen-Sznitman coalescent}
\author{Christina Goldschmidt and James B.\ Martin}
\date{}

\newcommand{\E}[1]{\ensuremath{\mathbb{E} \left[#1 \right]}}
\newcommand{\Prob}[1]{\ensuremath{\mathbb{P} \left(#1 \right)}}
\newcommand{\var}[1]{\ensuremath{\mathrm{var} \left(#1 \right)}}

\newcommand{\N}{\ensuremath{\mathbb{N}}}
\newcommand{\Ch}[2]{\ensuremath{\begin{pmatrix} #1 \\ #2 \end{pmatrix}}}
\newcommand{\ch}[2]{\ensuremath{\left( \begin{smallmatrix} #1 \\ #2
\end{smallmatrix} \right)}}

\newcommand{\fl}[1]{\ensuremath{\lfloor #1 \rfloor}}

\renewcommand{\subset}{\subseteq}
\newcommand{\convdistn}{\ensuremath{\stackrel{d}{\rightarrow}}}
\newcommand{\convprob}{\ensuremath{\stackrel{p}{\rightarrow}}}
\newcommand{\equidist}{\ensuremath{\stackrel{d}{=}}}
\newcommand{\bigO}{\ensuremath{\mathcal{O}}}
\newcommand{\argmax}{\ensuremath{\operatornamewithlimits{argmax}}}
\newcommand{\tE}{\ensuremath{\tilde{E}}}

\newtheorem{thm}{Theorem}[section]

\newtheorem{lemma}[thm]{Lemma}
\newtheorem{prop}[thm]{Proposition}
\newtheorem{cor}[thm]{Corollary}

\numberwithin{equation}{section}

\begin{document}
\maketitle

%%%%%%%%%%%%
% Abstract %
%%%%%%%%%%%%

\begin{abstract}
\noindent
We describe a representation of the Bolthausen-Sznitman coalescent in
terms of the cutting of random recursive trees.  Using this
representation, we prove results concerning the final collision of the
coalescent restricted to $[n]$: we show that the distribution of the
number of blocks involved in the final collision converges as
$n\to\infty$, and obtain a scaling law for the sizes of these blocks.
We also consider the discrete-time Markov chain giving the number of
blocks after each collision of the coalescent restricted to $[n]$; we
show that the transition probabilities of the time-reversal of this
Markov chain have limits as $n \rightarrow \infty$.  These results can
be interpreted as describing a ``post-gelation'' phase of the
Bolthausen-Sznitman coalescent, in which a giant cluster containing
almost all of the mass has already formed and the remaining small
blocks are being absorbed.
\end{abstract}

%%%%%%%%%%%%%%%%
% Introduction %
%%%%%%%%%%%%%%%%

\section{Introduction}

The Bolthausen-Sznitman coalescent, $\left(\Pi(t), t\geq 0\right)$, is a
Markov process which takes its values in the set of partitions of
$\N$.  It is most easily defined via its restriction
$\left(\Pi^{[n]}(t), t\geq 0\right)$ to the set $[n]:=\{1,2,\dots,
n\}$, for $n\geq 1$.  Denote by $\#\Pi^{[n]}(t)$ the number of blocks
of $\Pi^{[n]}(t)$. Then $\left(\Pi^{[n]}(t), t\geq 0\right)$ is a
continuous-time Markov chain whose transition rates are as follows: if
$\#\Pi^{[n]}(t)=b$, then any $k$ of the blocks present coalesce at
rate
\begin{equation} \label{eqn:lambdas}
\lambda_{b,k}=\frac{(k-2)!(b-k)!}{(b-1)!}, \qquad 2\leq k\leq b\leq n.
\end{equation}
It is usual to start the coalescent from the partition into
singletons, 
\[
\Pi(0)=\big(\{1\}, \{2\}, \{3\},\dots\big),
\]
and in this case we say that the coalescent is \textit{standard}.

The Bolthausen-Sznitman coalescent was first introduced in
\cite{Bolthausen/Sznitman}, in the context of the
Sherrington-Kirkpatrick model for spin glasses.  In
\cite{PitmanLambdaCoal}, Pitman demonstrated a great number of its
properties. He introduced the class of coalescents with multiple
collisions (also known as $\Lambda$-coalescents), gave a construction
of them based on Poisson random measures and studied the
Bolthausen-Sznitman coalescent as a member of this class.  Bertoin and
Le Gall \cite{Bertoin/LeGall} give an alternative derivation in terms
of the genealogy of a continuous-state branching process. Marchal
\cite{MarchalVienna} gives a construction via
regenerative sets.

In this paper, we describe a new representation for the
Bolthausen-Sznitman coalescent in terms of the cutting of random
recursive trees.

We then use this representation to prove results about the last
collision of the coalescent restricted to $[n]$, as $n \rightarrow
\infty$.  We obtain scaling laws for the sizes of the blocks involved
in the final collision; essentially, this collision involves one large
block and one or more smaller blocks, whose combined size behaves like
$n^U$, where $U$ has the uniform distribution on $[0,1]$.  We also
show that the distribution of the number of blocks involved in the
final collision converges as $n\to\infty$ (for example, the
probability that exactly two blocks are involved converges to $\log
2$).

More generally, we can also consider the discrete-time Markov chain
giving the number of blocks after each collision of the coalescent
restricted to $[n]$.  We show that the transition probabilities of the
time-reversal of this Markov chain have limits as $n \rightarrow
\infty$, which we make explicit.  (We observe in passing that the form
of these limiting probabilities yields certain infinite product
expansions of powers of $e$, which appear to be new).

These results can be interpreted as describing a ``post-gelation''
phase of the Bolthausen-Sznitman coalescent, in which a giant cluster
containing almost all of the mass has already formed and the very
small left-over blocks are being absorbed.  We also note that the tree
representation has an intrinsic asymmetry which contrasts strongly
with the exchangeability properties of the coalescent itself (for
example, the root of the tree always represents the block containing
$1$). This makes it possible to read off properties concerning a
tagged particle in the coalescent process directly from the tree
representation.

%%%%%%%%%%%%%%%%%%%%%%%%%%
% Random recursive trees %
%%%%%%%%%%%%%%%%%%%%%%%%%%

\section{Random recursive trees}

% ------------------ %
% The representation %
% ------------------ %

\subsection{The representation}
\label{sec:BS/RRT}

\begin{figure}[htb]
\begin{center}
\includegraphics{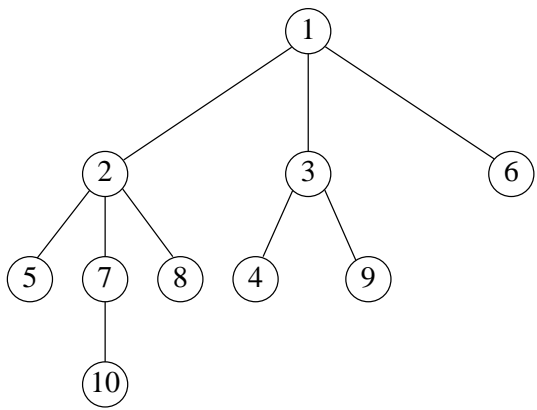}
\caption{A recursive tree on $[10]$.}
\label{fig:rrt}
\end{center}
\end{figure}

A tree on $n$ vertices labelled $1,2,\ldots,n$ is called a
\emph{recursive tree} if the vertex labelled $1$ is the root and, for
all $2 \leq k \leq n$, the sequence of vertex labels in the path from
the root to $k$ is increasing (Stanley~\cite{Stanley1} calls this an
unoriented increasing tree).  See Figure~\ref{fig:rrt} for an example
of a recursive tree.  Call a \emph{random recursive tree} a tree chosen
uniformly at random from the $(n-1)!$ possible recursive trees on $n$
vertices.  A random recursive tree can also be constructed as follows.
The vertex $1$ is distinguished as the root.  We imagine the vertices
arriving one by one.  For $k \geq 2$, vertex $k$ attaches itself to a
vertex chosen uniformly at random from $1, 2, \ldots, k-1$.  For a
detailed survey of results on recursive trees, see Smythe and
Mahmoud~\cite{Smythe/Mahmoud}.

We can represent an infinite tree by the infinite sequence of integers
$\{p_i\}_{i \geq 2}$, where $p_i$ is the parent of node $i$.  
(The condition for this to be a recursive tree is $1 \leq p_i < i$ for $i
\geq 2$.)

For the purposes of this article, it will be convenient also to define
a random recursive tree on a label set $\{l_1, l_2, \ldots, l_b\}$
where $l_1, l_2, \ldots, l_b$ are the blocks of a partition of $[n]$
for some $n$, listed in increasing order of least elements.  (That is,
we require that $l_1 < l_2 < \ldots < l_b$ where ``$<$'' is the total
order induced by the least elements of the blocks.)  The tree is
constructed in the obvious way: $l_1$ labels the root and $l_k$ is
attached to a vertex chosen uniformly at random from those labelled
$l_1, l_2, \ldots, l_{k-1}$.  Call the \emph{weight} of a label the
number of integers it contains and let $\mathcal{L}_n$ be the set of
partitions of $[n]$.

\begin{figure}
\begin{center}
\includegraphics{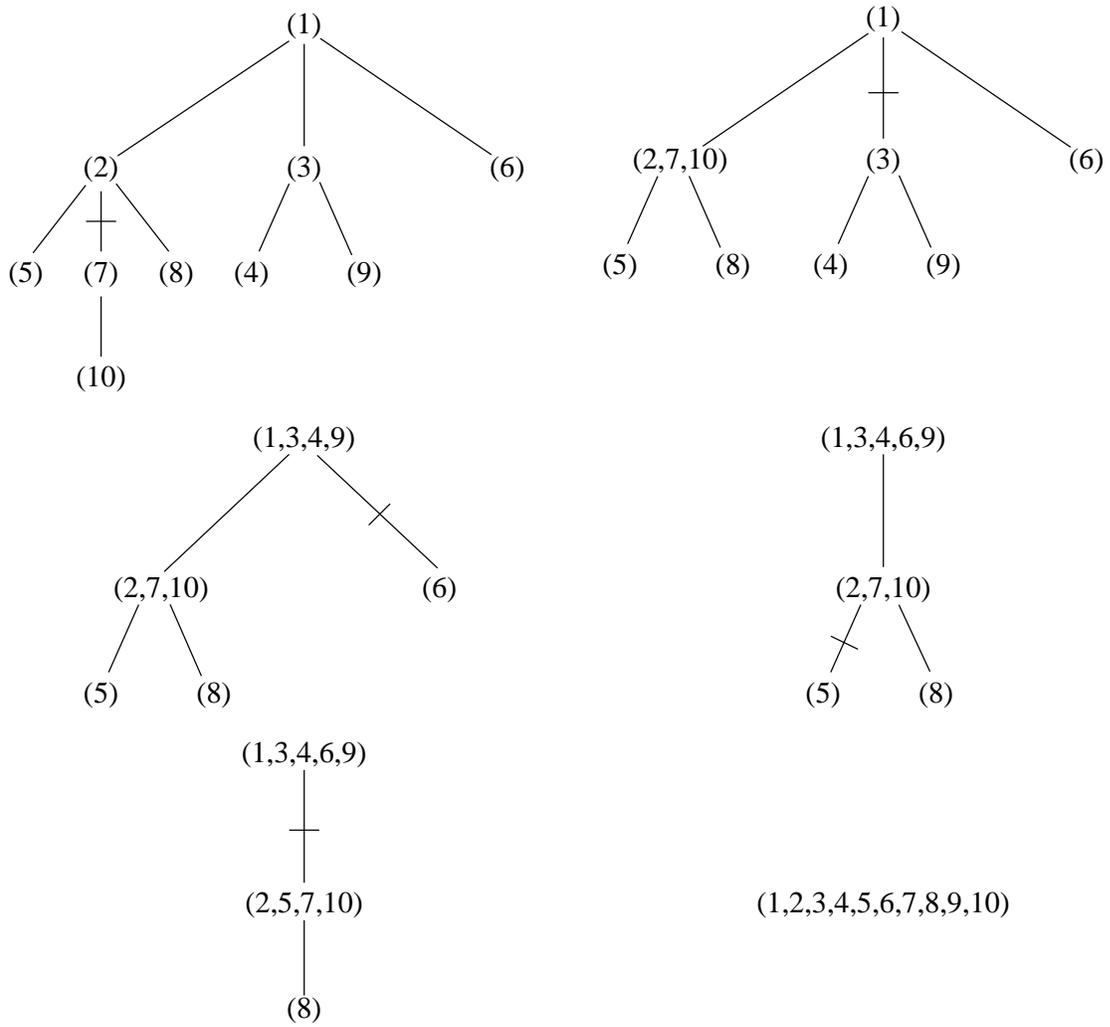}
\caption{The cutting procedure for a random recursive tree on $[10]$,
with the cuts indicated.}
\label{fig:rrtcut}
\end{center}
\end{figure}

Meir and Moon~\cite{Meir/Moon} define a cutting procedure which they
apply to random recursive trees.  They take a random recursive tree on
$[n]$ and pick an edge $e$ uniformly at random from the $n-1$ present.
This edge is deleted along with the entire subtree below it.  These
operations are then repeated until the root is isolated. The idea of
cutting combinatorial trees in this way appears to have been
introduced in Meir and Moon~\cite{Meir/MoonRandomTrees} and remains a
current research topic. Interest has tended to focus on the number of
cuts required to isolate the root in different types of trees.  Recent
references include Janson \cite{JansonPreprint,JansonVienna}, Fill,
Kapur and Panholzer \cite{Fill/Kapur/Panholzer} and Panholzer
\cite{PanholzerDMTCS,PanholzerVienna}.  In particular,
\cite{PanholzerVienna} treats the case of random recursive trees; the
author's presentation of this work at the MathInfo 2004 conference in
Vienna stimulated the present work.  A variant of the cutting
procedure will be the basis of our representation of the
Bolthausen-Sznitman coalescent.  Suppose that instead of throwing the
subtree below $e$ away, we add its labels to those of the vertex above
$e$.  We repeat this procedure until only the root remains, labelled
by $[n]$.  See Figure~\ref{fig:rrtcut} for an example.

\begin{prop} \label{prop:rrt}
Suppose $T$ is a random recursive tree on $L =
\{l_1,l_2,\ldots,l_b\} \in \mathcal{L}_n$.  Pick an edge at
random, cut it and add the labels below the cut to the label above.
Then the resulting tree is a random recursive tree on the new
label-set.
\end{prop}

\begin{proof}
(Essentially due to Meir and Moon~\cite{Meir/Moon}.)  The resulting
tree is clearly recursive because its labels still increase along all
paths away from the root.  So we need to show that it is chosen
uniformly from the set of all recursive trees with the same label-set.
Put another way, we will show that each of the $(b-1)! (b-1)$
recursive trees on the label-set $L$ with a single marked edge
corresponds to a tree constructed as follows: for some $2 \leq k \leq
b$, pick $k$ of the labels, say $l_{i_1}, l_{i_2}, \ldots, l_{i_k}$
(taken to be in increasing order), make a recursive tree on $L
\setminus \{l_{i_2}, \ldots, l_{i_k}\}$, make another recursive tree
on $\{l_{i_2}, \ldots, l_{i_k}\}$ and then join them together with an
edge between the vertices labelled $l_{i_1}$ and $l_{i_2}$.

There are $\ch{b}{k}$ ways of picking the $k$ labels $l_{i_1}, \ldots,
l_{i_k}$.  There are $(k-2)!$ ways of arranging the $k-1$ largest into
a recursive tree rooted at $l_{i_2}$.  There are $(b-k)!$ ways of
arranging the $b-k+1$ other labels into a recursive tree.  Clearly
each of the trees constructed in this way is distinct and also a
recursive tree.  The number which can be constructed is
\[
\sum_{k=2}^{b} \Ch{b}{k} (k-2)! (b-k)! = b! \sum_{k=2}^{b}
\frac{1}{k(k-1)} = (b-1)! (b-1).
\]
Hence, the claimed correspondence holds.
\end{proof}

\begin{prop}
Start with a random recursive tree on $[n]$ and associate an
independent exponential random variable with mean $1$ to each edge.
This exponential time is the time at which the edge is deleted, at
which point the labels in the subtree below it are instantaneously
added to the label of the vertex above the edge.  Then the set of
labels forms a partition of $[n]$ which evolves according to the
dynamics of the Bolthausen-Sznitman coalescent restricted to $[n]$.
\end{prop}

\begin{proof}
We need to show that the rate of coalescence of any set of $k$ of the
labels is $\lambda_{b,k}$ whenever there are $b$ vertices in the tree,
for $\lambda_{b,k}$ defined as at (\ref{eqn:lambdas}).  The total rate
of events when there are $b$ vertices is $b-1$.  The probability that
the next event will coalesce a particular $k$-set is worked out in the
same way as in the proof of Proposition~\ref{prop:rrt}.  Suppose we
start with label-set $L = \{l_1, l_2, \ldots, l_b\}$ and we want the
probability that the next event is the coalescence of $\{l_{i_1},
\ldots, l_{i_k}\}$.  There are $(k-2)!$ ways of making a recursive
tree on $\{l_{i_2}, \ldots, l_{i_k}\}$; there are $(b-k)!$ ways of
making a recursive tree on the remaining labels.  There are
$(b-1)!(b-1)$ recursive trees on a label-set of size $b$ with a single
marked edge and so the probability that we coalesce $\{l_{i_1},
\ldots, l_{i_k}\}$ is
\[
\frac{(k-2)! (b-k)!}{(b-1)!(b-1)}.
\]
Hence, the \emph{rate} at which we coalesce any $k$-set is
\[
\frac{(k-2)! (b-k)!}{(b-1)!}.
\]
The evolution is Markovian because, by Proposition~\ref{prop:rrt}, the
resulting tree is another random recursive tree, this time on $b-k+1$
labels.
\end{proof}

Because of the recursive way in which the original tree is built, the
representations are consistent for different $n$: if $\Pi^{[n]}(t)$ is
the partition given by the tree representation on $[n]$ at time $t$
then for all $t\geq0$ and all $n \in \N$,
\[
\Pi^{[n+1]}(t)|_{[n]}=\Pi^{[n]}(t).
\]
Thus we are able to define $(\Pi(t),t\geq0)$, the Bolthausen-Sznitman
coalescent on the whole of $\N$, by means of the cutting procedure
applied to an infinite random recursive tree (indexed by $\N$).

% ---------------------------------------------- %
% Random recursive trees, the Chinese restaurant %
% process and random permutations                %
% ---------------------------------------------- %

\subsection{Random recursive trees, the Chinese restaurant process and
random permutations}

There is a well-known correspondence between random recursive trees
and uniform random permutations.  We find it most convenient to
demonstrate this connection via the Chinese restaurant process of
Dubins and Pitman (see Chapter 3 of Pitman~\cite{PitmanStFl}), which
we now introduce.

For $0 \leq \alpha < 1$ and $\theta > -\alpha$, the Chinese restaurant
process with parameters $\alpha$ and $\theta$ (referred to hereafter
as the $(\alpha,\theta)$-Chinese restaurant process) is constructed as
follows. Person 1 enters a Chinese restaurant and sits at the first
table.  Person 2 sits either at the same table or at a new one; in
general, each subsequent person (numbered successively $3,4,\ldots,n$)
sits either at one of the occupied tables or at a new one.  Suppose
that people $1,2,\ldots,m$ have sat at $k$ tables where table $i$ has
$m_i$ customers (and $\sum_{i=1}^{k} m_i = m$).  Then person $m+1$
starts a new table with probability
\[
\frac{\theta + \alpha k}{m + \theta}
\]
and sits at table $i$ with probability
\[
\frac{m_i - \alpha}{m + \theta},
\]
for $1 \leq i \leq k$.  Once all $n$ customers have arrived, the
tables constitute the blocks of a partition of $[n]$ which are
consistent for different $n$.  Letting $n \rightarrow \infty$, these
blocks have \emph{asymptotic frequencies}, where the asymptotic
frequency of a block $B \subset \N$ is defined to be
\[
\lim_{n \rightarrow \infty} \frac{|B \cap [n]|}{n}.
\]
The existence of this limit for blocks generated by the Chinese
restaurant process is guaranteed by Kingman's theory of exchangeable
random partitions \cite{Kingman1, Kingman2}.  Moreover, if $(F_1, F_2,
\ldots)$ is the vector of asymptotic frequencies for a $(\alpha,
\theta)$-Chinese restaurant process in the order of the tables then
$(F_1, F_2, \ldots)$ has the Griffiths-Engen-McCloskey
$\mathrm{GEM}(\alpha,\theta)$ distribution.  If this vector is put into
decreasing order of size then it has the Poisson-Dirichlet
$\mathrm{PD}(\alpha,\theta)$ distribution.  (See Pitman~\cite{PitmanStFl}
or Arratia, Barbour and Tavar\'e~\cite{Arratia/Barbour/Tavare} for
more details.)

For $\alpha = 0$ and $\theta > 0$, there is a simpler construction.
Person 1 sits at table 1.  For $m \geq 1$, person $m+1$ starts a new
table with probability $\theta / (\theta + m)$ and sits to the left of
each of the $m$ existing customers with probability $1/(\theta + m)$.
In this case, the blocks created contain more information than is
needed to just give a partition of $[n]$.  In fact, they are the
cycles of a permutation of $[n]$, written in standard cycle notation
(i.e.\ the permutation maps $i$ to the customer seated to the left of
$i$).  The special case $\theta = 1$ gives a uniform random
permutation of $[n]$ (that is, a permutation chosen uniformly at
random from all the possible permutations of $[n]$).  Moreover, the
permutations are consistent for different $n$ and, when $n \rightarrow
\infty$, we obtain a uniform random permutation of $\N$.

We can find a $(0,1)$-Chinese restaurant process in the construction of
the random recursive tree on $[n+1]$.  For this purpose, it is easier
to imagine the random recursive tree labelled by the set
$\{0,1,\ldots,n\}$ rather than $\{1,2,\ldots,n+1\}$.  The root
(labelled 0) is fixed and does not appear in the permutation.
Vertices attached to the root correspond to individuals who start a
new table.  Thus, individual $1$ necessarily starts a new table.  If
vertex $k$ arrives and attaches to a vertex other than the root (say
$j$) then individual $k$ sits directly to the left of individual $j$.
As vertex $k$ is equally likely to attach to each of the vertices
labelled $0,1, \ldots, k-1$, individual $k$ is equally likely to sit
to the left of any of the individuals $1,2, \ldots, k-1$ or to form a
new table.

Thus, a random recursive tree on $[n+1]$ corresponds to a random
permutation of $[n]$.

Just to give some interpretation to the cutting procedure in the
Chinese restaurant, we embellish the usual model as follows.  Each
customer present in the restaurant has friends arrive at rate $1$ and
there is also a rate $1$ stream of customers who know no-one in the
restaurant at their time of arrival.  Customers who arrive and have a
friend present always sit to the left of that friend.  Friendless
customers sit at new tables.  Stop the process when there have been
$n$ arrivals.  A meal in the restaurant costs one euro.  At rate $1$,
each customer decides to leave and go home.  Whenever he leaves, any
of his friends who arrived after him want to leave (and their friends,
and so on).  (If the person who decided to depart was $k$ then it is
the subtree rooted at $k$ in the random recursive tree which departs.)
Anyone leaving gives the money for their meal to the person to whose
left $k$ sat down on entering the restaurant, so that he can pay for
them when he leaves.  If $k$ was, in fact, the first person to have
sat at that table then he collects all the money for the whole table
(plus, of course, the price of his own meal), takes it to the cashier
and departs.  Note that at any time $t$, the amount of money in the
cash register is the same as the weight of the label at the root in
the random recursive tree (i.e.\ the size of the block containing $1$
in the Bolthausen-Sznitman coalescent).  In this paper, we will answer
questions such as: how much does the last customer to leave pay the
cashier?

% ---------------------- %
% Using the construction %
% ---------------------- %

\subsection{Using the construction}

The construction of Subsection~\ref{sec:BS/RRT} gives an oddly
asymmetrical way of looking at the Bolthausen-Sznitman coalescent.
Rather than the usual exchangeable partitions approach here, because
the blocks are automatically ordered according to their least
elements, we have a size-biased view (smaller numbers are more likely
to be the smallest number in their blocks).  Moreover, a particular
realisation of a random recursive tree corresponds to a particular
conditioning of the coalescent (e.g.\ if $2$ and $5$ are both children
of $1$ then we are conditioning $2$ and $5$ only to be in the same
block when they have both coalesced with $1$).  Some facts are very
easy to read off from the tree and we will describe here some examples.

In Subsection 2.4 of \cite{PitmanLambdaCoal}, Pitman discusses the
size of the block containing $1$ at time $t$.  In terms of the random
recursive tree, in order to find the size of the block containing $1$
at time $t$, it suffices to know only the sizes of the descendencies
of all the children of the root, plus the times at which the edges
from the root to those children are severed.  By the connection with
the Chinese restaurant process, it is clear that in the limit as $n
\rightarrow \infty$ the vector of asymptotic frequencies of the
descendencies of the children of the root has the $\mathrm{GEM}(0,1)$
distribution (i.e.\ the $\mathrm{PD}(0,1)$ distribution with the
blocks in size-biased rather than decreasing order).  Moreover, by
construction, the times at which these descendencies are added to the
root are independent and identically distributed standard exponential
random variables.  Thus, we have here given a very straightforward
proof of Pitman's Corollary 16, reproduced below:

\begin{thm}[Pitman]
Let $\tilde{f}_1(t)$ be the frequency of the block containing $1$ at
time $t$ in a standard Bolthausen-Sznitman coalescent.  Then
\begin{description}
\item[(i)] the process $(\tilde{f}_1(t), t \geq 0)$ is Markovian, with
the same distribution as the process $(\gamma(1 - e^{-t}) / \gamma(1),
t \geq 0)$ where $(\gamma(s), s \geq 0)$ is a gamma process, with
stationary independent increments and $\Prob{\gamma(s) \in dx} =
\Gamma(s)^{-1} x^{s-1} e^{-x} dx, x > 0$.
\item[(ii)] The distribution of $\tilde{f}_1(t)$ is $\mathrm{Beta}(1 -
e^{-t}, e^{-t})$, and the process $(-\log(1 -  \tilde{f}_1(t)), t \geq
0)$ has non-stationary independent increments.
\item[(iii)] Let $J_1 \geq J_2 \geq \ldots$ be the ranked magnitudes
of jumps of the process $(\tilde{f}_1(t), t \geq 0)$, and let $T_i$ be
the time when the jump of magnitude $J_i$ occurs.  Then the
distribution of the sequence $(J_1, J_2, \ldots)$ is
$\mathrm{PD}(0,1)$, and this sequence is independent of the $T_i$,
which are independent with standard exponential distribution.
\end{description}
\end{thm}

((i), (ii) and (iii) are equivalent by properties of the Dirichlet
random measure; (iii) is immediate from the random recursive tree
approach.)

In Bolthausen and Sznitman~\cite{Bolthausen/Sznitman} and
Pitman~\cite{PitmanLambdaCoal}, it is shown that the marginal
distribution of the asymptotic frequencies of the blocks of the
coalescent (ranked in decreasing order), at fixed time $t > 0$, is
$\mathrm{PD}(e^{-t},0)$.  Jason Schweinsberg has observed that the
recursive tree provides a simple way of proving this fact.  We work
via the Chinese restaurant process.  Let $E_1, E_2, \ldots$ be
independent and identically distributed standard exponential random
variables.  Fix a time $t$ and imagine constructing the random
recursive tree in such a way that the vertex $i$ arrives with the edge
above it marked if $E_i < t$ and unmarked otherwise, for all $i \geq
2$.  (These marks are the cuts, but we will find it convenient here
not to collapse the tree at the cuts but rather just to keep track of
where they are.)  Then vertex $i$ has the smallest label in its block
if there are no marks in the path from the root to $i$.  Otherwise,
$i$ is in the same block as the closest vertex to $i$ in the path from
the root to $i$ which has no cuts above it.  Suppose now that vertices
$1,2,\ldots,n$ have arrived (possibly with marks) in the construction
of the random recursive tree.  Suppose also that they form $k$ blocks
(after cutting) of sizes $n_1, \ldots, n_k$ (where $\sum_{i=1}^{k} n_i
= n$).  Then vertex $n+1$ has $n$ possible places to attach.  It
creates a new block if (i) it attaches to the smallest element of one
of the $k$ blocks and (ii) it doesn't have a mark. This happens with
probability $ke^{-t}/n$. It adds to the $i$th block (of size $n_i$) if
it attaches to one of the $n_i - 1$ non-smallest elements of the block
or if it attaches to the smallest element \emph{and} arrives with a
mark.  This happens with probability $(n_i - e^{-t})/n$.  But these
are exactly the probabilities in the $(e^{-t},0)$-Chinese restaurant
process and so we can conclude that the full asymptotic frequencies in
decreasing order must have $\mathrm{PD}(e^{-t},0)$ distribution.

% ------------------------------------------------------------- %
% Translating results for trees into results for the coalescent %
% ------------------------------------------------------------- %

\subsection{Translating results for trees into results for the coalescent}

As mentioned above, Panholzer~\cite{PanholzerVienna} has studied the
number of cuts necessary to isolate the root in a random recursive
tree on $[n]$.  In the context of the Bolthausen-Sznitman coalescent,
this corresponds to the number of collision events that take place
until there is just a single block.  Let $J_n$ be this number of
collisions.

\begin{thm}[Panholzer]
$J_n$ has moments and centred moments as follows: for $k \in \N$, as
$n \rightarrow \infty$,
\begin{align*}
\E{J_n^k} & = \frac{n^k}{\log^k n} 
  + \left( H_k + k + \sum_{l=1}^{k} \Psi(l) \right) \frac{n^k}{\log^{k+1} n} 
  + \bigO \left( \frac{n^k}{\log^{k+2} n} \right) \\
\intertext{and}
\E{|J_n - \E{J_n}|^k} 
& = \left( (-1)^k + \sum_{l=0}^{k-1} \Ch{k-1}{l} (-1)^{k-l-1} \Psi(l+1) \right)
    \frac{n^k}{\log^{k+1} n} + \bigO \left( \frac{n^k}{\log^{k+2} n} \right)
\end{align*}
where $H_n = \sum_{i=1}^{n} \frac{1}{i}$ and $\Psi(x) = \frac{d}{dx}
\log \Gamma(x)$.  Thus,
\[
\frac{\log n}{n} J_n \convprob 1,
\]
as $n \rightarrow \infty$.
\end{thm}

Panholzer notes that it is not possible to obtain a limiting
distribution of a centred, scaled version of $J_n$ by the method
of moments.

%%%%%%%%%%%%%%%%%%%%%%
% The last collision %
%%%%%%%%%%%%%%%%%%%%%%

\section{The last collision of the Bolthausen-Sznitman coalescent}

% ------------ %
% Main results %
% ------------ %

\subsection{Main results}

Let $M_n$ be the sum of the sizes (or \emph{masses}) of the blocks not
containing the integer $1$ in the last collision of the
Bolthausen-Sznitman coalescent restricted to $[n]$, and
let $B_n$ be the number of blocks involved in the final collision.
For example, in Figure \ref{fig:rrtcut} we have $n = 10$, $M_n=5$ 
and $B_n=3$.

\begin{thm} \label{thm:main}
As $n \rightarrow \infty$,
\begin{equation} \label{eqn:yule}
\left(\frac{\log M_n}{\log n}, B_n \right) \convdistn (U, 1 + Y(UE)),
\end{equation}
where $(Y(t))_{t \geq 0}$ is a standard Yule process, $U$ is uniform
on $[0,1]$, $E$ is standard exponential and $(Y(t))_{t \geq 0}$, $U$
and $E$ are independent.
\end{thm}

The proof of this theorem is quite long and so we defer it to
Section~\ref{sec:mainpf}.  Theorem~\ref{thm:main} tells us that the
final collision of the coalescent restricted to $[n]$ is between a
very large block containing almost all of the mass and a collection of
very small blocks whose total mass is roughly $n^{U}$.  
We can make the distribution of the number
of the small blocks in the last collision explicit as follows.

\begin{prop} \label{prop:distn}
For $m \geq 1$, 
\begin{equation*}
\Prob{Y(UE) = m} = \frac{1}{m} \sum_{k=1}^{m} \Ch{m}{k} (-1)^{k+1}
\log (k+1).
\end{equation*}
This distribution has infinite mean.
\end{prop}

\begin{proof}
Let $q_m(t) = \Prob{Y(t) = m}$ for $m \geq 1$.  It is easily shown
(for example, from the forward equations) that
\begin{equation*}
q_m(t) = (1 - e^{-t})^{m-1} - (1 - e^{-t})^m
\end{equation*}
and so, by Proposition~\ref{prop:exponentialexpansion} (see
Section~\ref{subsec:recurrenceanalysis}), we have
\begin{equation*}
q_m(t) =
\begin{cases}
e^{-t} & \text{if $m = 1$} \\
\sum_{k=2}^{m} \ch{m-2}{k-2} (-1)^{k} (e^{-(k-1)t} - e^{-kt})
       & \text{if $m \geq 2$.}
\end{cases}
\end{equation*}
Thus,
\begin{align*}
\Prob{Y(UE) = m} & = \E{q_m(UE)} \\
& = \begin{cases}
\E{e^{-UE}} & \text{if $m = 1$} \\
\sum_{k=2}^{m} \ch{m-2}{k-2} (-1)^{k} (\E{e^{-(k-1)UE}} - \E{e^{-kUE}})
       & \text{if $m \geq 2$.}
\end{cases}
\end{align*}
For any $\theta > 0$,
\[
\E{e^{-\theta UE}} = \int_0^1 \int_0^{\infty} e^{-\theta x u} e^{-x}
dx du = \frac{1}{\theta} \log(\theta + 1)
\]
and so
\[
\Prob{Y(UE) = 1} = \log 2 
\]
and, for $m \geq 2$,
\begin{align*}
\Prob{Y(UE) = m} & = \sum_{k=2}^m \Ch{m-2}{k-2} (-1)^k
\left(\frac{1}{k-1} \log k - \frac{1}{k} \log(k+1) \right) \\ 
& = \frac{1}{m} \sum_{k=1}^m \Ch{m}{k} (-1)^{k+1} \log(k+1).
\end{align*}
Finally, as $\E{Y(t)} = e^t$,
\[
\E{Y(UE)} = \int_0^1 \int_0^{\infty} e^{ux} e^{-x} dx du = \int_0^1
\frac{1}{1-u} du = \infty.
\]
\end{proof}

For example, putting $m=1$ gives $\Prob{B_n=2}\to\log 2$ as 
$n\to\infty$. The joint convergence given in Theorem
\ref{thm:main} makes it possible to 
condition on the value of $Y(UE)$ to obtain results such as
\begin{cor}
Conditional on the event $\{B_n = 2\}$,
\[
\frac{\log M_n}{\log n} \convdistn V,
\]
where $V$ has density $\frac{1}{(1 + v) \log 2}$ on $[0,1]$.
\end{cor}

\emph{Remarks. (a)} We have already mentioned the connection between
random recursive trees and uniform random permutations.  In a
permutation, $M_n$ represents the length of a uniformly-chosen cycle.
In fact, considerably stronger results hold on the distribution of the
cycle lengths.  Let $K_n(t)$ be the number of cycles of length less
than or equal to $n^t$.  Then the functional central limit theorem
of DeLaurentis and Pittel~\cite{DeLaurentis/Pittel} states that
\begin{equation} \label{eqn:functionalCLT}
\left(\frac{K_n(t) - t \log n}{\sqrt{\log n}}, \quad 0 \leq t \leq 1 \right)
\end{equation}
(a random element of $D[0,1]$) converges weakly to Brownian motion on
$[0,1]$.  (This theorem was extended by Hansen~\cite{Hansen} to the
case of a permutation with distribution given by the Ewens sampling
formula i.e.\ a permutation generated by the $(0,\theta)$-Chinese
restaurant process.  An alternative proof of her more general theorem,
which is more in the style of this paper, was given by Donnelly, Kurtz
and Tavar\'e~\cite{Donnelly/Kurtz/Tavare}.  See also Feng and
Hoppe~\cite{Feng/Hoppe} for a path-level large deviations principle
for the Ewens sampling formula.) \newline \emph{(b)} Take any rooted
tree $T$, random or deterministic.  Suppose that each edge $e$ has a
random variable $\lambda_e$ attached to it, where the values
$\lambda_e$ are independent and identically distributed with a
continuous distribution.  Say that $\lambda_e$ is a record if it is
the largest value in the path from the root to $e$.
Janson~\cite{JansonPreprint} shows that the distribution of the number
of records in $T$ is the same as that of the number of random cuts
required to isolate the root.  To see this, each time cut the edge
with the largest $\lambda_e$ amongst those remaining.  Then by
symmetry, we always cut a uniformly-chosen edge, and so we get the
cutting procedure.  Moreover, $\lambda_e$ is a record if and only if
its edge is cut; thus, there must be the same number of records as of
cut edges.  In the records setting, $B_n$ has the distribution of the
size of the maximal subtree of a random recursive tree on $[n]$
containing the smallest record $\lambda_*$ and only edges $e$ with
$\lambda_e < \lambda_*$.

\begin{prop}
If $A_n$ is the time to absorption for the Bolthausen-Sznitman
coalescent on $[n]$, so that
\[
A_n = \inf \{ t \geq 0 : \#\Pi^{[n]}(t) = 1 \},
\]
then
\[
A_n - \log \log n \convdistn -\log E,
\]
where $E$ has a standard exponential distribution.
\end{prop}

\begin{proof}
This is a corollary of Lemma~\ref{lem:2} in Section~\ref{sec:mainpf}
(to connect the notation, we have $A_n = E_*$).
\end{proof}

This proposition is straightforward, but it nicely illustrates the
idea that the Bolthausen-Sznitman coalescent on $\N$ ``only just
fails'' to come down from infinity (see Pitman~\cite{PitmanLambdaCoal}
and Schweinsberg~\cite{Schweinsberg}).

% ------------------------------- %
% The coalescent in reversed time %
% ------------------------------- %

\subsection{The coalescent in reversed time}

We can regard Theorem~\ref{thm:main} as a statement about the first
event in a time-reversed version of the Bolthausen-Sznitman
coalescent.  More precisely, let $X = (X_0, X_1, \ldots )$ be the
discrete-time Markov chain describing the evolution of the number of
blocks in the coalescent restricted to $[n]$, so that $X_i$ is the number
of blocks after $i$ collisions, for $i \geq 0$.  (We can do this
because the Bolthausen-Sznitman coalescent is homogeneous, in the
sense that the sizes of the blocks do not influence the dynamics of
the process.)  $X$ has transition probabilities
\begin{equation} \label{eqn:fwdtransprobs}
p_{b,b-k+1} = \frac{b}{(b-1)} \frac{1}{k(k-1)}, \qquad 2 \leq k \leq b,
\end{equation}
and is clearly decreasing.  Started from $X_0 = n$, for any finite
$n$, the chain hits $1$ almost surely.  Let $\hat{X}$ be the
time-reversal of $X$, where $\hat{X}$ starts from 1.  We show that the
reversed-time transition probabilities
\begin{align*}
\hat{p}_{l,m} := & \lim_{n \rightarrow \infty} \Prob{\text{$X$ enters
$l$ from $m$} | \text{$X_0 = n$ and $X$ hits $l$}} \\
= & \lim_{n \rightarrow \infty} \Prob{\text{$\hat{X}$ jumps from $l$ to $m$} |
\text{$\hat{X}_0 = 1$ and $\hat{X}$ hits $l$ and $n$}}
\end{align*}
exist and can be made explicit.

\begin{thm} \label{thm:transprobs}
We have
\begin{equation} \label{eqn:transprobs}
\hat{p}_{l,m} = \begin{cases}
                \displaystyle{\frac{1}{m-1} 
                \sum_{k=1}^{m-1} \Ch{m-1}{k} (-1)^{k+1} \log(k+1)}
		& \text{if $m > l=1$} \\
		\displaystyle{\frac{m}{(l-1)(m-l)(m-l+1)} 
                \frac{\displaystyle{\sum_{j=1}^{m-1}} \Ch{m-1}{j} (-1)^{j+1} \log(j+1)}
                     {\displaystyle{\sum_{k=1}^{l-1}} \Ch{l-1}{k} (-1)^{k+1} \log(k+1)}} 
		& \text{if $m > l \geq 2$.}
		\end{cases}
\end{equation}
\end{thm}

The $l=1$ case of this theorem is a direct consequence of
Theorem~\ref{thm:main} and Proposition~\ref{prop:distn}.  In
principle, one could extend the probabilistic proof that we give of
this case to study individual cases of higher $l$; however, we have
not been able to find such a proof for the general case, and so the
self-contained proof given in Section~\ref{subsec:recurrenceanalysis}
is via the analysis of certain recurrence formulae.

\emph{Remark.} We note that $\hat{p}_{1,2} \neq \hat{p}_{2,3}$ and so
the time-reversal of the Bolthausen-Sznitman coalescent is not
self-similar in the sense of Bertoin~\cite{Bertoin}.  See
Pitman~\cite{PitmanLambdaCoal}, Basdevant~\cite{Basdevant} and
Marchal~\cite{MarchalVienna} for results on
time-reversing the Bolthausen-Sznitman coalescent to obtain a
fragmentation process.

% ------------------------ %
% Number theoretic results %
% ------------------------ %

\subsection{Number theoretic results}

Before proceeding to the proofs of Theorems~\ref{thm:main} and
\ref{thm:transprobs}, we observe in passing that
Theorem~\ref{thm:transprobs} entails various infinite product
expansions of powers of $e$, which seem not to have been previously
known.

\begin{thm} \label{thm:infiniteproducts}
For integers $r \geq 1$,
\begin{equation*}
e^{1/r} = \prod_{n=r}^{\infty} \left( \prod_{k=1}^{n} (k+1)^{\ch{n}{k}
(-1)^{k+1}} \right)^{1/(n-r+1)}.
\end{equation*}
For example, the $r=1$ case put more clearly reads
\begin{equation*}
e = \left( \frac{2}{1} \right)^{1/1} \left( \frac{2^2}{1 \cdot 3}
\right)^{1/2} \left( \frac{2^3 \cdot 4}{1 \cdot 3^3} \right)^{1/3}
\left( \frac{2^4 \cdot 4^4}{1 \cdot 3^6 \cdot 5} \right)^{1/4} \cdots
\end{equation*}
and the $r=2$ case reads
\begin{equation*}
e^{1/2} = \left( \frac{2^2}{1 \cdot 3}
\right)^{1/1} \left( \frac{2^3 \cdot 4}{1 \cdot 3^3} \right)^{1/2}
\left( \frac{2^4 \cdot 4^4}{1 \cdot 3^6 \cdot 5} \right)^{1/3}  \left(
\frac{2^5 \cdot 4^{10} \cdot 6}{1 \cdot 3^{10} \cdot 5^5} \right)^{1/4}
\cdots.
\end{equation*}
\end{thm}

The $r=1$ case of this theorem was first proved by Guillera and Sondow
(cited in Sondow~\cite{SondowFaster}) and bears a certain resemblance
to Pippenger's product \cite{Pippenger} for $e$.

\begin{proof}
Consider first the $r=1$ case.  In Theorem~\ref{thm:transprobs}, we
find the distribution $(\hat{p}_{1,n+1}, n \geq 1)$.  Summing over
$n$, we obtain
\begin{equation*}
\sum_{n=1}^{\infty} \frac{1}{n} \sum_{k=1}^{n} \Ch{n}{k} (-1)^{k+1}
\log(k+1) = 1.
\end{equation*}
Exponentiating both sides, we have 
\begin{equation*}
e = \prod_{n=1}^{\infty} \left( \prod_{k=1}^{n} (k+1)^{\ch{n}{k}
(-1)^{k+1}} \right)^{1/n},
\end{equation*}
which is Guillera and Sondow's product.

We now proceed by induction on $r$.  Assume the result up to $r-1$.
We have
\begin{equation*}
\sum_{m=r+1}^{\infty} \hat{p}_{r,m} = 1
\end{equation*}
and so, from (\ref{eqn:transprobs}),
\begin{equation*}
\sum_{m=r+1}^{\infty} \frac{m}{(r-1)(m-r)(m-r+1)} 
                \frac{\sum_{j=1}^{m-1} \ch{m-1}{j} (-1)^{j+1} \log(j+1)}
                     {\sum_{k=1}^{r-1} \ch{r-1}{k} (-1)^{k+1} \log(k+1)} 
= 1.
\end{equation*}
Hence,
\begin{align*}
& (r-1) \sum_{k=1}^{r-1} \Ch{r-1}{k} (-1)^{k+1} \log(k+1) \\
& \qquad = \sum_{m=r+1}^{\infty} \left( \frac{r}{m-r} - \frac{r-1}{m-r+1}
\right) \sum_{j=1}^{m-1} \Ch{m-1}{j} (-1)^{j+1} \log(j+1) \\
& \qquad = r \sum_{n=r}^{\infty} \frac{1}{n-r+1} \sum_{j=1}^{n}
\Ch{n}{j} (-1)^{j+1} \log(j+1) \\
& \qquad \qquad - (r-1) \sum_{n=r}^{\infty}
\frac{1}{n-r+2} \sum_{j=1}^{n} \Ch{n}{j} (-1)^{j+1} \log(j+1).
\end{align*}
Shifting the last term to the left-hand side and using the induction
hypothesis, we obtain
\begin{equation*}
1 = r \sum_{n=r}^{\infty} \frac{1}{n-r+1} \sum_{j=1}^{n}
\Ch{n}{j} (-1)^{j+1} \log(j+1)
\end{equation*}
and taking the exponential of both sides gives the result.
\end{proof}

% --------------------- %
% Proof of main theorem %
% --------------------- %

\subsection{Proof of Theorem~\ref{thm:main}} \label{sec:mainpf}

It is convenient to think of the formation of the random recursive
tree as occurring in continuous time, according to a linear birth
process with immigration.  (This is the approach used by Donnelly,
Kurtz and Tavar\'e~\cite{Donnelly/Kurtz/Tavare} when dealing with
random permutations.)  More specifically, we imagine that the root is
present at time 0 and that it has a special r\^ole.  Children of the
root ``immigrate'' at rate 1 and each immigrant initiates a family
which evolves as a standard Yule process.  Let $R_i$ be the time of
the $i$th immigration and let $P(t)$ be the number of immigrations up
to time $t$ (this is just a standard Poisson process).  Let
$Y_i(t)$ be the size of the the $i$th family at time $R_i + t$ and let
\[
T_n = \inf \left\{ t \geq 0: \sum_{i=1}^{P(t)} Y_i(t - R_i) = n-1 \right\},
\]
the time at which the total population present (including the root)
reaches $n$.  (Of course, $1 + \sum_{i=1}^{P(t)} Y_i(t - R_i)$ is
itself just a standard Yule process, but in what follows we find it
helpful to distinguish the root.)  Finally, let $C_n = P(T_n)$, the
number of immigrations up to (and including) time $T_n$.  Note that
$P$ and $Y_1, Y_2, \ldots$ are independent.

The theorem concerns what happens when we cut the tree so much that
all of the children of the root become severed.  We will think of the
cutting as occurring in continuous time (although subsequent to the
formation of the tree).  To this end, we associate an independent
standard exponential lifetime to each individual, which gives the time
at which the edge above it is cut. 

For each $i \geq 1$ we define a family of processes,
$\{(Y^{\text{thinned}(u)}_i(s), s \geq 0 ), u\geq0\}$, in the
following way.  The tree representing the $i$th Yule process at time
$s$ has size $Y_i(s)$. From this tree, we remove all the edges whose
cutting time is less than or equal to $u$ (and all the subtrees below
such edges). Let $Y^{\text{thinned}(u)}_i(s)$ be the size of the tree
that remains. Note that for any $u$, the process
$(Y_i^{\text{thinned}(u)}(s), s\geq 0)$ has the distribution of
$\left(Y(e^{-u}s), s\geq0\right)$ where $Y$ is a standard Yule
process.

For $1 \leq i \leq C_n$, let $E_i$
be the exponential lifetime associated with the $i$th child of the
root.  The last child to be cut is severed at time
\[
E_* = \max_{1 \leq i \leq C_n} E_i
\]
and is the $C_*$th child, where
\[
C_* = \argmax_{1 \leq i \leq C_n} E_i.
\]
$E_*$ is the time at which the root is isolated, and is therefore the
time of the last event in the coalescent process.  Let $R_* =
R_{C_*}$, the time at which the $C_*$th child arrived in the formation
of the tree.

In this framework, we have that 
\[
M_n = Y_{C_*}(T_n - R_*),
\]
the total family size (in the unthinned tree) of
individual $C_*$ (which is the last child of the root to be cut).
Similarly,
\[
B_n-1=Y^{\text{thinned}(E_*)}_{C_*}(T_n - R_*),
\]
the number of nodes still remaining in the (thinned) 
family of $C_*$ at the moment when it is cut. 

Let $t_n^{-} = \log n - (\log n)^{1/2}$ and $t_n^+ = \log n + (\log
n)^{1/2}$.  Let
\begin{gather*}
E_*^- = \max_{1 \leq i \leq P(t_n^-)} E_i, \\
C_*^- = \argmax_{1 \leq i \leq P(t_n^-)} E_i \\
\intertext{and}
R_*^- = R_{C_*^-}.
\end{gather*}
For completeness, if $P(t_n^-)=0$ (which becomes very unlikely as $n$
grows) then let $E_*^-=0$, $C_*^-=0$, and set $R_*^-$ equal to a
uniform random variable on $[0,t_n^-]$.

Now define
\begin{equation*}
\mathcal{G}_n = 
\left\{ T_n \in (t_n^-, t_n^+),\ C_* = C_*^-\
\text{ and }\
Y_{C_*^-}^{\text{thinned}(E_*^-)}(t_n^--R_*^-)
= Y_{C_*^-}^{\text{thinned}(E_*^-)}(t_n^+-R_*^-)
\right\}.
\end{equation*}
On the ``good set'' $\mathcal{G}_n$, we have $E_* = E_*^-$ and $R_* =
R_*^-$ almost surely.  Also, $R_*^- \sim \mathrm{U}[0, t_n^-]$ and is
independent of $E_*^-$ and of the Yule processes $Y_i$ for $i \geq 1$.
We will prove the following three lemmas:

\begin{lemma} \label{lem:1}
As $n \rightarrow \infty$,
\begin{gather*}
\frac{\log Y_{C_*^-}(t_n^- - R_*^-)}{t_n^- - R_*^-}\convdistn 1, \\ 
\frac{\log Y_{C_*^-}(t_n^+ - R_*^-)}{t_n^+ - R_*^-}\convdistn 1, \\
\intertext{and}
\frac{t_n^+ - R_*^-}{t_n^- - R_*^-}\convdistn 1.
\end{gather*}
\end{lemma}

\begin{lemma} \label{lem:2}
For any $n$,
\[
t_n^- e^{-E_*^-} \equidist E\wedge t_n^-,
\]
where $E$ is a standard exponential random variable.
\end{lemma}

\begin{lemma} \label{lem:3}
As $n \rightarrow \infty$, $\Prob{\mathcal{G}_n} \rightarrow 1$.
\end{lemma}

To show how the lemmas give the theorem, we proceed as follows.  On
$\mathcal{G}_n$, we have
\begin{gather} \label{eqn:cond1}
\frac{\log Y_{C_*^-}(t_n^- - R_*^-)}{\log n}\ 
\leq\ \frac{\log M_n}{\log n}\
\leq\ \frac{\log Y_{C_*^-}(t_n^+ - R_*^-)}{\log n} \\
\intertext{and}
\label{eqn:cond2}
B_n-1=Y_{C_*^-}^{\text{thinned}(E_*^-)}(t_n^- -R_*^-).
\end{gather}
So to show that $\left(\log M_n/\log n, B_n\right) \convdistn \left(U,
1+Y(UE)\right)$, it is enough to show that
\begin{equation}
\label{eqn:aimfor}
\left(
\frac{\log Y_{C_*^-}(t_n^- - R_*^-)}{\log n}, 
\frac{\log Y_{C_*^-}(t_n^+ - R_*^-)}{\log n},
Y_{C_*^-}^{\text{thinned}(E_*^-)}(t_n^- -R_*^-)
\right)
\convdistn \left(U,U,Y(UE)\right).
\end{equation}

We can write $(t_n^- -R_*^-)/t_n^- = U$ and $t_n^- e^{-E_*^-}= E\wedge
t_n^-$, where $U\sim U[0,1]$ and $E\sim \text{Exp}(1)$, and $U$, $E$
and $Y_1, Y_2, \ldots$ are independent.

Using Lemma \ref{lem:1} and the fact that $t_n^-/\log n \to 1$, we can
rewrite the left-hand side of (\ref{eqn:aimfor}) as
\[
\left( A_n U, B_n U, Y_{C_*^-}^{\text{thinned}(E_*^-)}(t_n^- U)
\right),
\]
where $A_n\convdistn 1$ and $B_n \convdistn 1$ as $n\to\infty$.
Thus, using Slutsky's Lemma, it is enough to show that
\begin{equation}\label{eqn:aimfor2}
\left(
U, U, 
Y_{C_*^-}^{\text{thinned}(E_*^-)}(t_n^- U)
\right)
\convdistn \left(U, U, Y(UE)\right).
\end{equation}
Recall that we can replace the process
$(Y_{C_*^-}^{\text{thinned}(E_*^-)}(s), s \geq 0)$ by
$(Y(e^{-E_*^-}s), s \geq 0)$, and that $t_n^- e^{-E_*^-}=
E\wedge t_n^-$. Then the left-hand side of (\ref{eqn:aimfor2}) becomes
\[
\left(U, U, Y\!\left(U (E\wedge t_n^-)\right) \right).
\]
Since $t_n^-\to\infty$, this indeed converges in distribution to
$\left(U,U,Y(UE)\right)$, as desired.

It remains to prove the three lemmas.

\emph{Proof of Lemma~\ref{lem:1}.} 
For a standard Yule process $Y$,
\[
\Prob{Y(t) = m} = (1 - e^{-t})^{m-1} - (1 - e^{-t})^{m}
\]
and so, for any $0 < \epsilon < 1$, 
\[
\Prob{ \frac{\log Y(t)}{t} \in (1 - \epsilon, 1 +
\epsilon]} = (1 - e^{-t})^{\fl{e^{(1 - \epsilon)t}}} - (1 -
e^{-t})^{\fl{e^{(1 + \epsilon)t}}},
\]
which converges to 1 as $t \rightarrow \infty$.  Thus,
\[
\frac{\log Y(t)}{t}\convdistn 1,
\]
as $t\to\infty$.  Hence, if $V_n$ is a sequence of random variables
such that $V_n\convdistn \infty$ as $n\to\infty$ (in the sense that
$\Prob{V_n>x}\to 1$ as $n\to\infty$, for all $x$), then also
\[
\frac{\log Y(V_n)}{V_n}\convdistn 1,
\]
as $n\to\infty$.

The first two parts of the lemma follow, since $R_*^-\sim U[0,t_n^-]$,
and so $t_n^- - R^-_*$ and $t_n^+ - R^-_*$ both converge to infinity
in distribution.

The last part is immediate from the facts that $R_*^-\sim U[0,t_n^-]$
and that $t_n^-/t_n^+\to 1$. \hfill $\Box$

\emph{Proof of Lemma~\ref{lem:2}.}  
Think of $P$ now as a marked Poisson process: points arrive at rate 1
and the $i$th point carries the mark $E_i$, where the $E_i$ are
independent and identically distributed with standard exponential
distribution.  So, for $s>0$, points with mark greater than or equal
to $s$ arrive as a Poisson process of rate $e^{-s}$.  The quantity
$E^-_*$ is the largest mark out of all those points arriving in the
interval $[0,t_n^-]$.  So
\begin{equation*}
\Prob{E_*^-<s} = \Prob{\text{no mark $\geq s$ arrives in $[0,t_n^-]$}}
= \exp(-t_n^-e^{-s}).
\end{equation*}
So for $x<t_n^-$, we have
\begin{equation*}
\Prob{t_n^- e^{-E_*^-} > x}
= \Prob{E_*^-<-\log \left(\frac{x}{t_n^-}\right)}
= e^{-x}.
\end{equation*}
Also, $E_*^-\geq 0$ with probability 1, so $t_n^- e^{-E_*^-}\leq t_n^-$
with probability 1.  Thus, indeed, $t_n^- e^{-E_*^-}$ has the
distribution of $E\wedge t_n^-$, as required. \hfill $\Box$

\emph{Proof of Lemma~\ref{lem:3}.}  The time, $T_n$, of the $(n-1)$th
birth in a standard Yule process satisfies
\[
T_n \equidist \tE_1 + \frac{1}{2} \tE_2 + \cdots + \frac{1}{n-1} \tE_{n-1},
\]
where $\tE_1, \tE_2, \ldots, \tE_{n-1}$ are independent and identically
distributed standard exponential random variables.  Hence,
$\E{T_n}=\log n + \bigO(1)$ and $\var{T_n} = \bigO(1)$, so Chebyshev's
inequality gives
\begin{equation} \label{eqn:T_nconv}
\Prob{T_n \in (t_n^-, t_n^+)} \to 1,
\end{equation}
as $n \to \infty$.  

Suppose now that $T_n \in (t_n^-, t_n^+)$ but that $C_*\neq C_*^-$.
This implies the event that, out of all the children of the root
arriving in $[0,t_n^+]$, the one with the largest exponential lifetime
arrived in $[t_n^-,t_n^+]$.  The probability of this event is simply
the ratio $(t_n^+ - t_n^-)/t_n^+$, which converges to 0 as
$n\to\infty$. Thus, we have
\begin{equation} \label{eqn:C*conv}
\Prob{C_* \neq C_*^-, T_n \in (t_n^-, t_n^+)} \to 0,
\end{equation}
as $n\to\infty$.

Finally, we need to show that 
\begin{equation}
\label{eqn:thirdneeded}
\Prob{Y_{C_*^-}^{\text{thinned}(E_*^-)}(t_n^--R_*^-)
= Y_{C_*^-}^{\text{thinned}(E_*^-)}(t_n^+-R_*^-)} \to 1.
\end{equation}
As observed earlier, we can replace the process
$Y_{C_*^-}^{\text{thinned}(u)}(s)$ by $Y(e^{-u}s)$ and so we need to
show that
\[
\Prob{Y(e^{-E_*^-}(t_n^- - R_*^-)) = Y(e^{-E_*^-}(t_n^+
- R_*^-))} \to 1.
\]
As before, we can put $t_n^- - R_*^- = t_n^- U$ and
$e^{-E_*^-}=(E\wedge t_n^-)/t_n^-$ to rewrite this as
\[
\Prob{Y\left(U (E\wedge t_n^-)\right)
= Y\left(U (E\wedge t_n^-) \tfrac{t_n^+ -R_*^-}{t_n^- -R_*^-}\right)}
\to 1.
\]
Now $\Prob{E > t_n^-}\to 0$ and $\big(t_n^+ -R_*^-\big)/\big(t_n^-
-R_*^-\big)\convdistn 1$ as $n\to\infty$ (by Lemma \ref{lem:1}), so it
will be enough to show that
\[
\Prob{\text{a jump of $Y$ occurs in the interval $(UE, (1+\epsilon)UE)$}}
\to 0, 
\]
as $\epsilon\downarrow 0$.  This is easily seen to be true (for
example, by integrating over $U$ and $E$), since $\E{Y(t)}$ is
continuous in $t$.

Putting (\ref{eqn:T_nconv}), (\ref{eqn:C*conv}) and
(\ref{eqn:thirdneeded}) together, we have that $\Prob{\mathcal{G}_n}
\to 1$ as $n \to \infty$, as required.~$\hfill\Box$

% ----------------------------------------- %
% Proof of transition probabilities theorem %
% ----------------------------------------- %

\subsection{Proof of Theorem~\ref{thm:transprobs}}
\label{subsec:recurrenceanalysis}

The self-contained proof of Theorem~\ref{thm:transprobs} which we give
here mostly consists of the analysis of certain recurrence formulas
via generating functions (see Wilf~\cite{Wilf} for an excellent
introduction to this method).  We will start by proving the $l=1$ case
and will then see that this enables us to prove the rest of the
theorem.  Let
\begin{equation*}
x_n^{(m)} = \Prob{\text{$X$ enters $1$ from $m$} | X_0 = n}
\end{equation*}
for $m \geq 2$ (we do not need to condition on $X$ entering $1$ ever
as this occurs almost surely).  Then,
\begin{equation} \label{eqn:initialcond}
x_m^{(m)} = p_{m, 1} = \frac{1}{(m-1)^2}
\end{equation}
and, by the Markov property, for $n \geq m+1$,
\begin{equation} \label{eqn:recurrence}
x_n^{(m)} = \sum_{k=2}^{n-m+1} p_{n, n-k+1}\ x_{n-k+1}^{(m)} 
            = \frac{n}{n-1} \sum_{k=2}^{n-m+1}
            \frac{x_{n-k+1}^{(m)}}{k(k-1)}.
\end{equation}
We can re-phrase the $l = 1$ case of Theorem~\ref{thm:transprobs} as

\begin{thm} \label{thm:l=1}
Let $m \geq 2$.  If $(x_n^{(m)}, n \geq m)$ satisfies
(\ref{eqn:initialcond}) and (\ref{eqn:recurrence}) then
\begin{equation*}
\lim_{n \rightarrow \infty} x_n^{(m)} = \frac{1}{m-1}
\sum_{k=1}^{m-1} \Ch{m-1}{k} (-1)^{k+1} \log(k+1).
\end{equation*}
\end{thm}

The generating function of the sequence $(x_n^{(m)}, n \geq m)$ we
will use is defined as
\begin{equation} \label{eqn:genfun}
f_{(m)}(t) = \sum_{n=1}^{\infty} x_{n+m-1}^{(m)} t^{n+m-2},
\end{equation}
for $t \in [0,1)$.

The key to proving Theorem~\ref{thm:l=1} is the following result on
the convergence of power series.

\begin{prop} \label{prop:Littlewood}
Suppose that $f(t) = \sum_{n=1}^{\infty} f_n t^n$ is a power series
with positive real coefficients such that 
\begin{equation} \label{eqn:incrementbound}
|f_n - f_{n-1}| \leq \frac{C}{n-1}
\end{equation}
for some constant $C \in (0, \infty)$ and all $n \geq 2$.  Suppose
that 
\begin{equation}
(1 - t) f(t) \rightarrow f^{*}
\end{equation}
as $t \rightarrow 1-$.  Then $\lim_{n \rightarrow \infty} f_n = f^{*}$.
\end{prop}

\begin{proof}
We have $(1-t)f(t) = \sum_{n=1}^{\infty} (f_n - f_{n-1})t^n$, where we
take $f_0 = 0$ for convenience.  Then by Littlewood's Tauberian
theorem (see, for example, Theorem 1.1 of Chapter 2 in van de
Lune~\cite{vandeLune})
\begin{equation*}
f^{*} = \sum_{n=1}^{\infty} (f_n - f_{n-1}).
\end{equation*}
This sum telescopes and so we have $\lim_{n \rightarrow \infty} f_n =
f^{*}$, as required.
\end{proof}

In view of Proposition~\ref{prop:Littlewood}, we would like to prove a
condition like (\ref{eqn:incrementbound}) for the sequence
$(x^{(m)}_n, n \geq m)$.  It turns out to be easier to prove this for
a more general situation.  Let the sequence $(x_n, n \geq m)$ be
defined as follows for fixed $m \geq 2$: for some $x_m \in (0,1]$,
\begin{equation} \label{eqn:generalrecurrence}
x_n = \sum_{k=2}^{n-m+1} \frac{a_{n,k}}{b_n} x_{n-k+1}, \quad  n \geq
m+1,
\end{equation}
where the non-negative coefficients $(a_{n,k})_{2 \leq k \leq n}$
satisfy
\begin{equation} \label{eqn:consistencyconds}
a_{n,k} = \tfrac{n-k+1}{n+1} a_{n+1,k} + \tfrac{k+1}{n+1} a_{n+1,k+1},
\qquad 2 \leq k \leq n
\end{equation}
and
\begin{equation} \label{eqn:b_ns}
b_n = a_{n,2} + a_{n,3} + \cdots + a_{n,n}, \qquad n \geq 2.
\end{equation}
Conditions (\ref{eqn:consistencyconds}) and (\ref{eqn:b_ns}) are
exactly those satisfied by the rates in a general $\Lambda$-coalescent
and (\ref{eqn:generalrecurrence}) is the recurrence that arises for
its reversed-time transition probabilities.  In Pitman's notation of
\cite{PitmanLambdaCoal}, $a_{n,k} = \ch{n}{k} \lambda_{n,k}$ and $b_n
= \sum_{k=2}^{n} \ch{n}{k} \lambda_{n,k}$, where $\lambda_{n,k} =
\int_0^1 x^{k-2} (1 - x)^{n-k} \Lambda(dx)$ and $\Lambda$ is any
finite measure.  In the Bolthausen-Sznitman case, $a_{n,k} =
\frac{n}{k(k-1)}$ and $b_n = n-1$.

\begin{lemma} \label{lem:b_ns}
For $n \geq 3$, $b_n = b_{n-1} + \tfrac{2}{n} a_{n,2}$.
\end{lemma}

\begin{proof}
We have
\begin{align*}
b_{n-1} 
& = a_{n-1,2} + a_{n-1,3} + \cdots + a_{n-1,n-1} \\
& = \left( \tfrac{n-2}{n} a_{n,2} + \tfrac{3}{n} a_{n,3} \right) +
    \left( \tfrac{n-3}{n} a_{n,3} + \tfrac{4}{n} a_{n,4} \right) \\
& \qquad + \cdots + \left( \tfrac{2}{n} a_{n,n-2} + \tfrac{n-1}{n}
  a_{n,n-1} \right) + \left( \tfrac{1}{n} a_{n,n-1} + a_{n,n} \right)
  \qquad \text{by (\ref{eqn:consistencyconds})} \\
& = -\tfrac{2}{n} a_{n,2} + a_{n,2} + a_{n,3} + \cdots + a_{n,n} \\
& = b_n - \tfrac{2}{n} a_{n,2}.
\end{align*}
\end{proof}

\begin{lemma} \label{lem:increment}
For $n \geq m+2$,
\begin{equation} \label{eqn:increment}
x_{n-1} - x_n = \frac{1}{b_n} \left\{ \sum_{k=2}^{n-m} \tfrac{n-k}{n}
a_{n,k} (x_{n-k} - x_{n-k+1}) - \tfrac{m-1}{n} a_{n,n-m+1} x_m \right\}. 
\end{equation}
\end{lemma}

\begin{proof}
We have
\begin{align*}
b_n x_n 
& = a_{n,2} x_{n-1} + a_{n,3} x_{n-2} + \cdots + a_{n,n-m+1} x_m \\
\intertext{and}
b_{n-1} x_{n-1} 
& = a_{n-1,2} x_{n-2} + a_{n-1,3} x_{n-3} + \cdots +
    a_{n-1,n-m} x_{m} \\
& = \left( \tfrac{n-2}{n} a_{n,2} + \tfrac{3}{n} a_{n,3} \right) x_{n-2} +
    \left( \tfrac{n-3}{n} a_{n,3} + \tfrac{4}{n} a_{n,4} \right) x_{n-3} \\
& \qquad + \cdots + \left( \tfrac{m}{n} a_{n,n-m} +
  \tfrac{n-m+1}{n}a_{n,n-m+1} \right) x_m \\
& = b_n x_n - a_{n,2} x_{n-1} - a_{n,3} x_{n-2} - \cdots 
    - a_{n,n-m+1} x_m \\
& \qquad 
  + \left( \tfrac{n-2}{n} a_{n,2} + \tfrac{3}{n} a_{n,3} \right) x_{n-2} +
  \cdots + \left( \tfrac{m}{n} a_{n,n-m} +
  \tfrac{n-m+1}{n}a_{n,n-m+1} \right) x_m \\
& = b_n x_n - \tfrac{2}{n} a_{n,2} x_{n-1} + \tfrac{n-2}{n} a_{n,2}
(x_{n-2} - x_{n-1}) + \tfrac{n-3}{n} a_{n,3} (x_{n-3} - x_{n-2}) \\
& \qquad + \cdots + \tfrac{m}{n} a_{n,n-m} (x_m - x_{m+1}) - \tfrac{m-1}{n}
a_{n,n-m+1} x_m.
\end{align*}
Hence, 
\begin{align*}
& \left( b_{n-1} + \tfrac{2}{n} a_{n,2} \right) x_{n-1} \\
& \qquad = b_n x_n + \tfrac{n-2}{n} a_{n,2} (x_{n-2} - x_{n-1}) 
  + \cdots + \tfrac{m}{n} a_{n,n-m} (x_m - x_{m+1}) 
  - \tfrac{m-1}{n} a_{n,n-m+1} x_m
\end{align*}
and so, by Lemma~\ref{lem:b_ns},
\begin{equation*}
x_{n-1} - x_{n} 
= \frac{1}{b_n} \left\{ \tfrac{n-2}{n} a_{n,2}
  (x_{n-2} - x_{n-1}) + \cdots + \tfrac{m}{n} a_{n,n-m} (x_m - x_{m+1}) 
  - \tfrac{m-1}{n} a_{n,n-m+1} x_m \right\}
\end{equation*}
which gives (\ref{eqn:increment}).
\end{proof}

\begin{lemma} \label{lem:1/nbound}
For $n \geq m+1$, we have
\begin{equation} \label{eqn:1/nbound}
\left| x_{n-1} - x_{n} \right| \leq \frac{m}{n-1}.
\end{equation}
\end{lemma}

\begin{proof}
We proceed by induction.  For $n = m+1$, by
(\ref{eqn:generalrecurrence}) we have
\[
x_{m} - x_{m+1} = x_{m} \left(1 - \frac{a_{m+1,2}}{b_{m+1}} \right)
\leq 1.
\]
Now suppose that (\ref{eqn:1/nbound}) holds up to $n-1$.  From
(\ref{eqn:increment}) and this induction hypothesis we have
\begin{align*}
\left| x_{n-1} - x_n \right| 
& \leq \frac{1}{b_n} \left\{ \sum_{k=2}^{n-m} \tfrac{n-k}{n} a_{n,k}
\tfrac{m}{n-k} + \tfrac{m-1}{n} a_{n,n-m+1} \right\} \\
& = \frac{1}{b_n} \left\{ \tfrac{m}{n} (a_{n,2} + \cdots +
a_{n,n-m}) + \tfrac{m-1}{n} a_{n,n-m+1} \right\} \\
& \leq \frac{1}{b_n} \frac{m}{n} \sum_{k=2}^{n} a_{n,k} \\
& \leq \frac{m}{n-1}.
\end{align*}
Hence result.
\end{proof}

Thus, we know that for the sequences $(x_n^{(m)}, n \geq m)$
associated with the Bolthausen-Sznitman coalescent,
\[
|x^{(m)}_{n-1} - x^{(m)}_n| \leq \frac{m}{n-1}
\]
for all $n \geq m+1$.

As a final preparation for the proof of Theorem~\ref{thm:l=1}, we
mention a simple identity.

\begin{prop} \label{prop:exponentialexpansion}
For $m \geq 2$,
\begin{equation}
(1 - e^{-r})^{m-1} - (1 - e^{-r})^{m} = \sum_{k=2}^{m} \Ch{m-2}{k-2}
(-1)^{k} (e^{-(k-1)r} - e^{-kr}).
\end{equation}
\end{prop}

\begin{proof}
By induction on $m$.
\end{proof}

We are now ready to prove Theorem~\ref{thm:l=1}.

\emph{Proof of Theorem~\ref{thm:l=1}.}
Recall that for $m \geq 2$ and $t \in [0,1)$,
\begin{equation*}
f_{(m)}(t) = \sum_{n=1}^{\infty} x_{n+m-1}^{(m)} t^{n+m-2}.
\end{equation*}
Then
\begin{equation*}
f'_{(m)}(t) = \sum_{n=1}^{\infty} (n+m-2) x^{(m)}_{n+m-1} t^{n+m-3}.
\end{equation*}
From (\ref{eqn:recurrence}), we have 
\[
(n-1) x_n^{(m)} = n \sum_{k=2}^{n-m+1} \frac{x^{(m)}_{n-k+1}}{k(k-1)}
\]
and so
\begin{align}
t f'_{(m)}(t) 
& = \sum_{n=1}^{\infty} (n+m-2) x^{(m)}_{n+m-1} t^{n+m-2} \notag \\
& = \frac{(m-1) t^{m-1}}{(m-1)^2}
    + \sum_{n=2}^{\infty} (n+m-1) t^{n+m-2} 
         \sum_{k=2}^{n} \frac{x_{n+m-k}^{(m)}}{k(k-1)} \notag \\
& = \frac{t^{m-1}}{m-1}
    + \sum_{n=2}^{\infty} \sum_{k=2}^{n} 
      (n+m-k-1) x^{(m)}_{n+m-k}t^{n+m-k-1} 
      \cdot \frac{t^{k-1}}{k(k-1)} \notag \\
& \qquad + \sum_{n=2}^{\infty} \sum_{k=2}^{n} 
      x^{(m)}_{n+m-k} t^{n+m-k-1} \cdot \frac{t^{k-1}}{k-1}. \label{eqn:endup}
\end{align}
Now,
\begin{align*}
\sum_{k=1}^{\infty} \frac{t^{k}}{k(k+1)} & = \frac{1}{t} \left(
(1-t)\log(1-t) + t \right) \\
\intertext{and}
\sum_{k=1}^{\infty} \frac{t^k}{k} & = -\log(1 - t)
\end{align*}
and so (\ref{eqn:endup}) becomes
\[
t f'_{(m)}(t) = \frac{t^{m-1}}{m-1} + \left((1-t)\log(1-t) + t \right)
f'_{(m)}(t) - f_{(m)}(t) \log(1-t).
\]
Re-arranging and integrating gives
\begin{align*}
(1-t)f_{(m)}(t) 
& = - \frac{1}{m-1}\int_0^t \frac{u^{m-1}}{\log(1-u)} du \\
& = \frac{1}{m-1} \int_0^{-\log(1-t)} 
    \frac{(1 - e^{-r})^{m-1} - (1 - e^{-r})^{m}}{r} dr,
\end{align*}
where we have changed variable with $u = 1 - e^{-r}$.  By
Proposition~\ref{prop:exponentialexpansion}, we have
\begin{align*}
(1-t)f_{(m)}(t) 
& = \frac{1}{m-1} \sum_{k=2}^{m} \Ch{m-2}{k-2} (-1)^{k} 
    \int_0^{-\log(1-t)} \frac{(e^{-(k-1)r} - e^{-kr})}{r} dr.
\end{align*}
Letting $t \rightarrow 1-$, the integral on the right-hand side
becomes Frullani's integral and so
\begin{align}
\lim_{t \rightarrow 1-} (1 - t) f_{(m)}(t)
& = \frac{1}{m-1} \sum_{k=2}^{m} \Ch{m-2}{k-2} (-1)^k 
    \log \left(\frac{k}{k-1} \right) \notag \\
& = \frac{1}{m-1} \sum_{k=1}^{m-1} \Ch{m-1}{k} (-1)^{k+1}
    \log(k+1). \label{eqn:almostthere}
\end{align}
Thus, by Proposition~\ref{prop:Littlewood}, we have
\[
\hat{p}_{1,m} = \lim_{n \rightarrow \infty} x_n^{(m)} 
= \frac{1}{m-1} \sum_{k=1}^{m-1} \Ch{m-1}{k} (-1)^{k+1} \log(k+1).
\]
\hfill $\Box$

It remains only to show how the rest of Theorem~\ref{thm:transprobs} follows
from Theorem~\ref{thm:l=1}.  We first note that
\[
\Prob{\text{$X$ enters $l$ from $m$} | \text{$X_0 = n$ and $X$ hits $l$}}
= \frac{\Prob{\text{$X$ enters $l$ from $m$} | X_0 = n}}
       {\Prob{\text{$X$ hits $l$}| X_0 = n}}.
\]
Let
\[
y_n^{(l,m)} = \Prob{\text{$X$ enters $l$ from $m$} | X_0 = n}
\]
and
\[
y_n^{(l)} = \Prob{\text{$X$ hits $l$}| X_0 = n}.
\]
Then 
\[
y_m^{(l,m)} = \frac{m}{(m-1)} \frac{1}{(m-l+1)(m-l)}
\]
and the Markov property gives us that
\[
y_n^{(l,m)} = \frac{n}{n-1} \sum_{k=2}^{n-m+1}
\frac{y_{n-k+1}^{(l,m)}}{k(k-1)}, \qquad n \geq m+1.
\]
Thus, using (\ref{eqn:initialcond}) and (\ref{eqn:recurrence}),
$y_n^{(l,m)} = \frac{m(m-1)}{(m-l+1)(m-l)} x_n^{(m)}$ and so by
Theorem~\ref{thm:l=1},
\[
\lim_{n \rightarrow \infty} y_n^{(l,m)} = \frac{m}{(m-l+1)(m-l)}
\sum_{k=1}^{m-1} \Ch{m-1}{k} (-1)^{k+1} \log(k+1).
\]
Likewise, $y_l^{(l)} = 1$ and for $n \geq l+1$,
\[
y_n^{(l)} = \frac{n}{n-1} \sum_{k=2}^{n-l+1}
\frac{y_{n-k+1}^{(l)}}{k(k-1)}.
\]
Thus, $y_n^{(l)} = (l-1)^2 x_n^{(l)}$ and so
\[
\lim_{n \rightarrow \infty} y_n^{(l)} = (l-1)
\sum_{k=1}^{l-1} \Ch{l-1}{k} (-1)^{k+1} \log(k+1).
\] 
Hence, finally,
\[
\hat{p}_{l,m} 
= \lim_{n \rightarrow \infty} \frac{y^{(l,m)}_n}{y^{(l)}_n} 
= \frac{m}{(l-1)(m-l)(m-l+1)} 
\frac{\sum_{j=1}^{m-1} \ch{m-1}{j} (-1)^{j+1} \log(j+1)}
     {\sum_{k=1}^{l-1} \ch{l-1}{k} (-1)^{k+1} \log(k+1)}.
\]

It is clear from our argument that $\hat{p}_{l,m}$ must be
non-negative for all $m > l \geq 1$.  However, by Fatou's lemma, we
know only that
\[
\sum_{m=l+1}^{\infty} \hat{p}_{l,m} \leq 1.
\]
We would like to show that $(\hat{p}_{l,m}, m \geq l+1)$ is really a
distribution (i.e.\ that this inequality is in fact an equality).  The
only thing that can go wrong is if mass ``escapes to infinity''.  To
show that this does not happen, we firstly make clear the probability
space on which we are working.  For each state $n$, sample $d_n$
independently from the distribution $(p_{n,k}, 1 \leq k \leq n-1)$.
Then for any $n$, the path of $X$ started with $X_0 = n$ is $n, d_n,
d_{d_n}, \ldots$.  Thus, we have coupled the paths of the Markov chain
$X$ started from all different possible states.  Now for fixed $l \geq
1$ consider the events
\[
E_n := \{ \text{$X$ enters $l$ from $n$} \} = \{d_n = l\},
\]
for $n \geq l+1$.  We have $\Prob{E_n} = p_{n,l} =
\frac{n}{(n-1)}\frac{1}{(n-l+1)(n-l)}$.  As $\sum_{n=l+1}^{\infty}
\frac{n}{(n-1)}\frac{1}{(n-l+1)(n-l)} < \infty$, by the Borel-Cantelli
Lemma, $\Prob{\text{$E_n$ i.o.}} = 0$.  Hence, we must have
\[
\sum_{m=l+1}^{\infty} \hat{p}_{l,m} = 1.
\]

\emph{Remark.} It does not appear that the methods of this section can
be extended to a more general $\Lambda$-coalescent, despite the
encouraging fact that Lemma~\ref{lem:1/nbound} holds for all
$\Lambda$-coalescents.  The generating function methods used rely
crucially on special structure of the Bolthausen-Sznitman case, in
particular the fact that $p_{n,n-k+1}$ decomposes neatly into two
simple factors, one involving only $n$ and the other involving only
$k$.

%%%%%%%%%%%%%%%%%%%%
% Acknowledgements %
%%%%%%%%%%%%%%%%%%%%

\section*{Acknowledgments}

The argument in the proof of Theorem~\ref{thm:l=1} is due to Noam
Elkies.  We should also like to thank Nevin Kapur for showing us an
equivalent argument. Many thanks to Philippe Flajolet for a helpful
discussion and to Jason Schweinsberg for his comments and suggestions.
C.G.\ would like to express her gratitude to Jean Bertoin for having
welcomed her to Paris 6 for a year and for his support, then and
subsequently.

%%%%%%%%%%%%%%
% References %
%%%%%%%%%%%%%%

%\bibliography{lambdacoalescent,recursivetrees}
\bibliography{recursivetrees}
\begin{tabbing}
\hspace{8cm} \= \kill
Christina Goldschmidt,      \> James B. Martin, \\
Pembroke College and Statistical Laboratory, 
                            \> LIAFA, \\ 
University of Cambridge,    \> CNRS et Universit\a'e Paris 7, \\
Centre for Mathematical Sciences,
                            \> 2 place Jussieu (case 7014), \\
Wilberforce Road,           \> 75251 Paris Cedex 05, \\
Cambridge                   \> France \\
CB3 0WB                     \> \\
UK                          \> \\
\texttt{C.Goldschmidt@statslab.cam.ac.uk} 
                            \> \texttt{James.Martin@liafa.jussieu.fr} \\
\texttt{http://www.statslab.cam.ac.uk/\~{}cag27}
                            \> \texttt{http://www.liafa.jussieu.fr/\~{}martin} 
\end{tabbing}

\end{document}